\documentclass[a4paper,11pt]{article}
\usepackage{fullpage}
\usepackage{graphics}
\usepackage[latin2]{inputenc}
\usepackage{multicol}
\usepackage{amsmath,amssymb,amsbsy,amsthm,amsfonts}
\usepackage{epsfig} 

\begin{document}
\def\note#1{\marginpar{\small #1}}

\def\tens#1{\pmb{\mathsf{#1}}}
\def\vec#1{\boldsymbol{#1}}

\def\norm#1{\left|\!\left| #1 \right|\!\right|}
\def\fnorm#1{|\!| #1 |\!|}
\def\abs#1{\left| #1 \right|}
\def\ti{\text{I}}
\def\tii{\text{I\!I}}
\def\tiii{\text{I\!I\!I}}

\def\diver{\mathop{\mathrm{div}}\nolimits}
\def\grad{\mathop{\mathrm{grad}}\nolimits}
\def\Div{\mathop{\mathrm{Div}}\nolimits}
\def\Grad{\mathop{\mathrm{Grad}}\nolimits}

\def\tr{\mathop{\mathrm{tr}}\nolimits}
\def\cof{\mathop{\mathrm{cof}}\nolimits}
\def\det{\mathop{\mathrm{det}}\nolimits}

\def\lin{\mathop{\mathrm{span}}\nolimits}
\def\pr{\noindent \textbf{Proof: }}
\def\pp#1#2{\frac{\partial #1}{\partial #2}}
\def\dd#1#2{\frac{\d #1}{\d #2}}

\def\T{\mathcal{T}}
\def\R{\mathbb{R}}
\def\bx{\vec{x}}
\def\be{\vec{e}}
\def\bef{\vec{f}}
\def\bec{\vec{c}}
\def\bs{\vec{s}}
\def\ba{\vec{a}}
\def\bn{\vec{n}}
\def\bphi{\vec{\varphi}}
\def\btau{\vec{\tau}}
\def\bc{\vec{c}}
\def\bg{\vec{g}}

\def\bW{\vec{W}}
\def\bT{\tens{T}}
\def\bD{\tens{D}}
\def\bF{\tens{F}}
\def\bB{\tens{B}}
\def\bV{\tens{V}}
\def\bS{\tens{S}}
\def\bI{\tens{I}}
\def\bi{\vec{i}}
\def\bv{\vec{v}}
\def\bfi{\vec{\varphi}}
\def\bk{\vec{k}}
\def\b0{\vec{0}}
\def\bom{\vec{\omega}}
\def\bw{\vec{w}}
\def\p{\pi}
\def\bu{\vec{u}}

\def\ID{\mathcal{I}_{\bD}}
\def\IP{\mathcal{I}_{p}}
\def\Pn{(\mathcal{P})}
\def\Pe{(\mathcal{P}^{\eta})}
\def\Pee{(\mathcal{P}^{\varepsilon, \eta})}

\def\Ln#1{L^{#1}_{\bn}}

\def\Wn#1{W^{1,#1}_{\bn}}

\def\Lnd#1{L^{#1}_{\bn, \diver}}

\def\Wnd#1{W^{1,#1}_{\bn, \diver}}

\def\Wndm#1{W^{-1,#1}_{\bn, \diver}}

\def\Wnm#1{W^{-1,#1}_{\bn}}

\def\Lb#1{L^{#1}(\partial \Omega)}

\def\Lnt#1{L^{#1}_{\bn, \btau}}

\def\Wnt#1{W^{1,#1}_{\bn, \btau}}

\def\Lnd#1{L^{#1}_{\bn, \btau, \diver}}

\def\Wntd#1{W^{1,#1}_{\bn, \btau, \diver}}

\def\Wntdm#1{W^{-1,#1}_{\bn,\btau, \diver}}

\def\Wntm#1{W^{-1,#1}_{\bn, \btau}}

\newtheorem{Theorem}{Theorem}[section]
\newtheorem{Example}{Example}[section]
\newtheorem{Lem}{Lemma}[section]
\newtheorem{Rem}{Remark}[section]
\newtheorem{Def}{Definition}[section]
\newtheorem{Col}{Corollary}[section]
\newtheorem{Proposition}{Proposition}[section]

\newcommand{\Om}{\Omega}
\newcommand{ \vit}{\hbox{\bf u}}
\newcommand{ \Vit}{\hbox{\bf U}}
\newcommand{ \vitm}{\hbox{\bf w}}
\newcommand{ \ra}{\hbox{\bf r}}
\newcommand{ \vittest }{\hbox{\bf v}}
\newcommand{ \wit}{\hbox{\bf w}}
\newcommand{ \fin}{\hfill $\square$}

\newcommand{\ZZ}{\mathbb{Z}}
\newcommand{\CC}{\mathbb{C}}
\newcommand{\NN}{\mathbb{N}}
\newcommand{\V}{\zeta}
\newcommand{\RR}{\mathbb{R}}
\newcommand{\EE}{\varepsilon}
\newcommand{\Lip}{\textnormal{Lip}}
\newcommand{\XX}{X_{t,|\textnormal{D}|}}
\newcommand{\PP}{\mathfrak{p}}
\newcommand{\VV}{\bar{v}_{\nu}}
\newcommand{\QQ}{\mathbb{Q}}
\newcommand{\HH}{\ell}
\newcommand{\MM}{\mathfrak{m}}
\newcommand{\rr}{\mathcal{R}}
\newcommand{\tore}{\mathbb{T}_3}
\newcommand{\Z}{\mathbb{Z}}
\newcommand{\N}{\mathbb{N}}

\newcommand{\F}{\overline{\boldsymbol{\tau} }}

\newcommand{\moy} {\overline {\vit} }
\newcommand{\moys} {\overline {u} }
\newcommand{\mmoy} {\overline {\wit} }
\newcommand{\g} {\nabla }
\newcommand{\G} {\Gamma }
\newcommand{\x} {{\bf x}}
\newcommand{\E} {\varepsilon}
\newcommand{\BEQ} {\begin{equation} }
\newcommand{\EEQ} {\end{equation} }
\makeatletter
\@addtoreset{equation}{section}
\renewcommand{\theequation}{\arabic{section}.\arabic{equation}}

\newcommand{\hs}{{\rm I} \! {\rm H}_s}
\newcommand{\esp} [1] { {\bf I} \! {\bf H}_{#1} }

\newcommand{\vect}[1] { \overrightarrow #1}

\newcommand{\hsd}{{\rm I} \! {\rm H}_{s+2}}

\newcommand{\HS}{{\bf I} \! {\bf H}_s}
\newcommand{\HSD}{{\bf I} \! {\bf H}_{s+2}}

\newcommand{\hh}{{\rm I} \! {\rm H}}
\newcommand{\lp}{{\rm I} \! {\rm L}_p}
\newcommand{\leb}{{\rm I} \! {\rm L}}
\newcommand{\lprime}{{\rm I} \! {\rm L}_{p'}}
\newcommand{\ldeux}{{\rm I} \! {\rm L}_2}
\newcommand{\lun}{{\rm I} \! {\rm L}_1}
\newcommand{\linf}{{\rm I} \! {\rm L}_\infty}
\newcommand{\expk}{e^{ {\rm i} \, {\bf k} \cdot \x}}
\newcommand{\proj}{{\rm I}Ê\! {\rm P}}

\renewcommand{\theenumi}{\Roman{section}.\arabic{enumi}}

\newcounter{taskcounter}[section]

\newcommand{\bib}[1]{\refstepcounter{taskcounter} {\begin{tabular}{ l p{13,5cm}} \hskip -0,2cm [\Roman{section}.\roman{taskcounter}] & {#1}
\end{tabular}}}

\renewcommand{\thetaskcounter}{\Roman{section}.\roman{taskcounter}}

\newcounter{technique}[section]

\renewcommand{\thetechnique}{\roman{section}.\roman{technique}}

\newcommand{\tech}[1]{\refstepcounter{technique} {({\roman{section}.\roman {technique}}) {\rm  #1}}}

\newcommand{\B}{\mathcal{B}}

\newcommand{\diameter}{\operatorname{diameter}}


%
%
%
%


%
%

\begin{center}
\Large{\textbf{ Theory for the Rotational Deconvolution model of Turbulence with Fractional regularization}}
\end{center}
\begin{center}
\textbf{{Hani Ali}}$^{\hbox{a}}$\\
$^{\hbox{a}}$\textit{Univ Paris-Sud, Laboratoire de Math\'ematiques d'Orsay,\\
 Orsay Cedex, F-91405\\
 Hani.Ali@math.u-psud.fr}
\end{center}
\bigskip
\textbf{Abstract}
\bigskip\\
We introduce a new regularization of the rotational Navier-Stokes equations that we call the Rotational Approximate Deconvolution Model (RADM). We generalize the deconvolution type model,  studied by Berselli and Lewandowski \cite{bresslilewandowski}, to the  RADM model with fractional regularization where
the convergence of the solution is studied with weaker conditions on the parameter regularization. 


\smallskip
\smallskip



\section{{Introduction}}

 \hskip 0.5cm Let us consider the rotational Navier-Stokes equations in a  three dimensional torus $\mathbb{T}_3$,
\begin{align}
\nabla \cdot \bv &=0, \label{nsBM}\\
\bv_{,t} -{\vec{v}^{}} \times  
 \nabla \times {\vec{v}^{}} -  \nu \Delta \bv  + \nabla P &=  {\bef},\label{nsBLM}
\end{align}
with initial data  $\bv(\bx, 0) = \bv_{0}(\bx)$. Here, $\bv$ is the fluid velocity field,  $P$ is the Bernoulli scalar  or the dynamic pressure $P=p+\frac{1}{2}\vec{u}^2$,
$\bef$ is the  external body forces and  $\nu$ represents  the
viscosity.

When the flow is turbulent it is well known that equations (\ref{nsBM})-(\ref{nsBLM})
are unstable in numerical simulations. Therefore,  the numerical turbulent models are needed for real
simulations of turbulent flows. 
The  alpha models is an example among  these numerical models that have  attracted a lot of interest from both the pure and applied mathematicians.  These models are based on a  filtering or averaging obtained
through the application of the inverse of the Helmholtz operator 
\begin{align}
\label{exact}
\mathbb{A} := I- \alpha^2 \Delta. \end{align}
Similar to the alpha models
\cite{CFHOTW99b, FDT02, ILT05, CHOT05, CLT06, Re08, LST11, ZhouFan1, ZhouFan2}, Olson and Titi \cite{OT2007} used a more general filtering  with fractional regularization given by 
\begin{align}
\label{frac}
 \bv= \mathbb{A}_{\theta}\overline{\bv}:= \left( I + \alpha^{2\theta}(-\Delta)^{\theta} \right) \overline{\bv}  \textrm { for some fixed  }  \theta >0,
\end{align}  
such that  $\bv$ denotes the unfiltered velocity and $\overline{\bv}$ denotes the filtered one.
This filter   is a differential filter \cite{germano}, that commutes with differentiation under periodic boundary conditions. 

The Approximate Deconvolution Model of turbulence (ADM) belongs also to the class of numerical turbulent models. A complete study of the ADM models has been carried out during the last years thanks to their well-posedness   and their good performance in numerical simulation \cite{adamsstolz, AdamStolz, AdamStolzKleiser}. The ADM model consists in
replacing the nonlinear  term in the Navier-Stokes equations $ \nabla \cdot ( \bv \otimes \bv ) $  by $  \overline{\nabla \cdot ( D_N \bv \otimes D_N \bv) }$, where the exact Helmholtz operator (\ref{exact}) is used to define  the deconvolution operator $D_N$ of order $N$  as follows:
\begin{align}
D_{N} = \sum_{i = 0}^N ( I - \mathbb{A}_{}^{-1} )^i.
\end{align}
 Existence, regularity and
uniqueness of a solution to this model for a general deconvolution of order $N$, were proved in \cite{dunca06}. The case $N = 0$ was studied  in details in \cite{LL03, LL06b}. Recently, Berselli and Lewandowski \cite{bresslilewandowski} studied  a generalized  ADM model that consists in replacing the deconvolution operator $D_N$ by a more general   deconvolution operator $D_{N,\theta}$, with fractional regularization. This deconvolution  operator  is defined by 
\begin{align}
D_{N,\theta} = \sum_{i = 0}^N ( I - \mathbb{A}_{\theta}^{-1} )^i.
\end{align}
For $\theta > \frac{3}{4}$, Berselli and Lewandowski  \cite{bresslilewandowski} showed  the convergence  of the ADM solution,  with fractional regularization, to a
solution of the average  Navier-Stokes Equations  when N goes to infinity. 
The existence and
the uniqueness of a solution to this model in  the case $N = 0$ and $\theta \ge \frac{1}{6}$ are studied  in  \cite{HA11b}.

In this paper, we introduce   a new model called the Rotational 
Approximate 
Deconvolution model of turbulence (RADM). 
 Let $\overline{\bv}$ be the filtered or the averaged velocity field which is given by the relation (\ref{frac}). In order  to obtain the RADM, we first write the Reynolds averaged rotational Navier-Stokes equations 
 
\begin{align}
\begin{split}
\overline{\bv}_{,t} - \overline{\bv \times \nabla \times  {\bv}} -  \nu \Delta \overline{\bv} +\nabla \overline{P}
 = \overline{{\bef}},\label{NOTCLOSEDBLM}
\end{split}
\end{align}
which are not in closed form due to the averaged nonlinear term  $\overline{\bv \times \nabla \times  {\bv}}$ that  we cannot write in terms of $\overline{\bv}$ alone.

As the  main goal of turbulence modeling is
to derive simplified and computationally realizable closure models,   we use the
deconvolution operator $D_{N,\theta}$ to approximate the velocity  $\bv$  with $ D_{N,\theta} \overline{\bv} $ and write the equation for $ \overline{\bv}$ as bellow  

\begin{align}
\begin{split}
\overline{\bv}_{,t} - \overline{D_{N,\theta} \overline{\bv} \times \nabla \times  D_{N,\theta} \overline{\bv}} -  \nu \Delta \overline{\bv} +\nabla \overline{P}
 = \overline{{\bef}}.\label{CLOSEDBLM}
\end{split}
\end{align}

This yields the initial value problem
\begin{align}
\nabla \cdot \bw &=0, \label{llBM}\\
\bw_{,t} - \overline{D_{N,\theta} \bw \times \nabla \times D_{N,\theta}\bw} -  \nu \Delta \bw + \nabla q
 &=  \overline{\bef},\label{llBLM}
\end{align}
with initial data $\bw(\bx, 0) = \bw_{0}(\bx)=\overline{\bv_{0}}$, considered in $(0,T)\times \mathbb{T}_3$ and subject to   periodic boundary conditions with mean value equal to zero.
Note that  the  couple  $(\bw,q)$ is  the approximation of  $ (\overline{\bv}, \overline{P})$.

When $N = 0$, the operator $D_{0,\theta} = I$.   Thus the zeroth order RADM is given by the following equations 
\begin{align}
\nabla \cdot \bw &=0, \label{zerollBM}\\
\bw_{,t} -  \overline{ \bw \times \nabla \times \bw} -  \nu \Delta \bw + \nabla q
 &=  \overline{\bef}.\label{zerollBLM}
\end{align}

The model (\ref{llBM})-(\ref{llBLM}) is a deconvolution model for the rotational form of the nonlinearity. Therefore
calling it the RADM model seems natural. This model is interesting as it can be discretized using methods similar to those used for the rotational form of the Navier-Stokes equations. Recent work on the discretizations of the  Navier-Stokes equations  \cite{Benzi07anefficient} revealed  that solving the NSE in
this form can lead to better efficiency, especially in the case of low viscosity. 



The purpose of this paper is to study the  RADM model with fractional regularization $\theta$ and look for the limiting case where we can
prove global existence and uniqueness of regular solutions $(\vec{w}_N, q_N )$ for a fixed $N$.
 The analysis of  the zeroth order RADM (\ref{zerollBM})-(\ref{zerollBLM}) is the same as the general case. For this reason, the regularization parameter $\theta$ does not depends on $ N$.  
Another goal of this paper is to  study  the question of the convergence of the RADM  solution  when $N$ goes to infinity. 
Namely,  we show that when N tends to infinity and $\frac{1}{6} \le \theta < 1$   
 the unique solution 
$(\vec{w}_N, q_N ) $ of the RADM model satisfies  
\begin{align}
\nabla \cdot  \bw_{} &=0, \label{intintAthetallBM}\\
\bw_{,t} -  \overline{\mathbb{A}_{\theta} \bw \times \nabla \times \mathbb{A}_{\theta}\bw} -  \nu \Delta \bw + \nabla q
 &=  \overline{\bef},\\
\bw(\bx, 0) = \bw_{0}(\bx)& =\overline{\bv_{0}}. \label{intintAthetallBLM}
\end{align}

This result  of consistency and stability   illustrates   how the system
approximates the rotational Navier-Stokes equations. 

  
The rest of the paper is organized as follows. In section 2, we start by giving some mathematical
tools such as the space of functions we are using and  the deconvolution operator $D_{N,\theta}$.  Section 3 provides  existence and uniqueness results in the case where $  0 \le N < \infty$.     The last section provides the main result of this paper which concerns the study of the convergence of the RADM solution   when $N$ goes to infinity. 
Note that the  results of this paper  require   weaker  conditions on the parameter $\theta$  than that required in \cite{bresslilewandowski}.

\section{Mathematical setting}
 \subsection{Notations}
 
\ In this section  we fix notation of function spaces that we shall  employ.\\ 
We denote by $L^p(\tore)$ and $H^{s}(\tore)$, $s \ge -1, \ 1 \le p \le \infty,$ the usual Lebesgue and Sobolev spaces over $\tore$, and define the Bochner spaces $C(0,T;X), L^p(0,T;X)$  in the standard way. 
The Sobolev spaces $ \vec{H}^{s}=H^s(\tore)^3$,  of mean-free functions are classically characterized in terms of
the Fourier series
$$\vec{H}^{s} = \left\{ \bv(\bx)=\sum_{{\vec{k} \in {\mathcal T}_{3}}}^{}\bc_{\vec{k}} e^{i\vec{k}\cdot \vec{x}}, \left(\bc_{\vec{k}} \right)^{*}=\bc_{-\vec{k}}, \bc_0 = 0,   \|\bv \|_{s,2}^2= \sum_{{\vec{k} \in {\mathcal T}_{3}}}^{}| \vec{k} |^{2s} |\bc_{\vec{k}}|^2<\infty  \right\},$$
where  $ \left(\bc_{\vec{k}}^N \right)^{*}$ denote the complex conjugate $\bc_{\vec{k}}^N.  $
In addition we introduce
\begin{align*}
\vec{H}^{s}_{\sigma }&=\left\{ \bv \in \vec{H}^s; \; \nabla \cdot \bv
=0 \textrm{ in }  \tore \right\},\\
\vec{H}^{-s}&=\left(\vec{H}^{s}_{} \right)^{'},\quad
\vec{L}^2_{}=\vec{H}^{0},\quad 
\vec{L}^2_{\sigma }=\vec{H}^{0}_{\sigma }.
\end{align*}
Throughout we will use $C$ to denote an arbitrary constant which may change from line to
line.

\subsection{The generalized deconvolution operator}\label{sec:dec}
The deconvolution operator considered in this paper is a generalized deconvolution operator that is  constructed by using  the Helmholtz equation with fractional regularization  \cite{bresslilewandowski,A01}.
Let $\alpha >0$, $s \ge -1$, $  0 \le \theta  \le 1$, $\bv \in \vec{H}^{s}$ and let  $\overline{\bv} \in \vec{H}^{s+2\theta}$ be the unique solution to the
 equations (compare to \cite{OT2007})
 \begin{align}
  \alpha^{2 \theta}(-\Delta)^{\theta} \overline{\bv} + \overline{\bv} = \bv, \label {VPYRA}\\
\nabla \cdot\vec{v}= \nabla \cdot\overline{\vec{v}}=0.
\end{align}
We also shall  denote by $\mathbb{A}_{\theta}$ the operator 
\BEQ \mathbb{A}_{\theta} : \begin{array} {l} \vec{H}^{s+2\theta} \longrightarrow {\vec{H}^{s}}, \\
\bv \longrightarrow   \alpha^{2 \theta}(-\Delta)^{\theta} \bv  + \bv. \end{array}  \EEQ

The non-local operator $\mathbb{A}_{\theta}$ is
defined through the Fourier transform
\begin{equation}
\widehat{\mathbb{A}_{\theta}{\bv(\bk)}}=\left(1+\alpha^{2\theta}|{\bk}|^{2\theta}\right)\widehat{\bv}({\bk}).
\end{equation}
Therefore, one has 
\BEQ \overline{\bv} = \mathbb{A}_{\theta}^{-1} \bv, \EEQ and   

\BEQ
\label{lemme2.2}
\| \overline{\bv}\|_{s+2\theta} \le \frac{1}{\alpha^{2\theta}} \|\bv \|_{s}.
\EEQ


Let us consider the operators 
$$ D_{N,\theta} = \sum_{i = 0}^N ( I - \mathbb{A}_{\theta}^{-1} )^i.$$

 When $\theta =1$ and for a fixed $N >0$, we recover the van Cittert deconvolution operator   used in  Large eddy simulation (LES)  
by Stolz and Adams \cite{adamsstolz} see also \cite{dunca06}.  This deconvolution operator, for different values of $\theta$,  was used in \cite{DuncaLewandowski} to study  the  rate of convergence of the ADM model 
to the mean Navier-Stokes Equations. 

A straightforward calculation yields

\BEQ D_{N,\theta}  \left ( \sum_{ {\bk}\in {\cal  I}_3} \bc_{\bk} \expk \right ) =  
\sum_{ {\bk}\in {\cal  I}_3}\left(1+ \alpha^{2 \theta} | {\bk} |^{2 \theta}\right) \left ( 1- \left ( \frac{\alpha^{2 \theta} | {\bk} |^{2 \theta}}
 {1 +  \alpha^{2 \theta} | {\bk} |^{2 \theta}}\right )^{N+1} \right )  \bc_{\bk} \expk. \EEQ 
Thus 
 \BEQ D_{N,\theta}  \left ( \sum_{ {\bk}\in {\cal  I}_3} \bc_{\bk} \expk \right ) =  \sum_{ {\bk}\in {\cal  I}_3}  \widehat{D_{N,\theta}}(\bk)  \bc_{\bk} \expk, \EEQ 
 where we have for all ${\bk}\in {\cal  I}_3,$
 
 \begin{align}
 \widehat{D_{0,\theta}}(\bk) &= 1, \\
 1 \le  \widehat{D_{N,\theta}}(\bk) &\le N+1 \quad \hbox { for each } N >0, \\
 \hbox{ and  }   \widehat{D_{N,\theta}}(\bk) &\le \widehat{\mathbb{A}_{\theta}} := \left(1+ \alpha^{2 \theta} | {\bk} |^{2 \theta}\right) \hbox{ for  a fixed } \alpha >0. \label{haniali}
  \end{align}
%
%
One can prove the following Lemma  (see in \cite{bresslilewandowski} ): 
\begin{Lem}
\label{lemme1}
For all  $s \ge -1$, ${\bk}\in {\cal  I}_3,$ and for each $N >0$ there exist a constant $C >0$ such that for all $\bv$ sufficiently smooth  we have 
\begin{align}
  &\|\bv\|_{s,2} \le    \|D_{N,\theta} \left ( \bv \right )\|_{s,2} \le  (N+1)  \|\bv\|_{s,2},\label{lemme13a}\\
    &\|\bv\|_{s,2} \le   C \|D_{N,\theta} \left ( \bv \right )\|_{s,2} \le  C \|{\mathbb{A}_{\theta}}^{\frac{1}{2}}D_{N,\theta}^{\frac{1}{2}}(\bv)\|_{s,2},\label{lemme13b}\\
     &\|{\mathbb{A}_{\theta}}^{\frac{1}{2}}D_{N,\theta}^{\frac{1}{2}}(\overline{\bv})\|_{s,2}  \le  \|\bv\|_{s,2}, \label{lemme13}\\    
      &\|\bv\|_{s+\theta,2} \le   C(\alpha) \|{\mathbb{A}_{\theta}}^{\frac{1}{2}}D_{N,\theta}^{\frac{1}{2}}(\bv)\|_{s,2}.\label{lemme13d}
  \end{align}
\end{Lem}

\section{The RADM model: existence and uniqueness of a  weak solution} 

Let us consider the following  Rotational Approximate Deconvolution Model (RADM)

\begin{align}
\nabla \cdot  \bw_{} &=0, \label{thetallBM}\\
\bw_{,t} - \overline{D_{N,\theta} \bw \times \nabla \times D_{N,\theta}\bw} -  \nu \Delta \bw + \nabla q
 &=  \overline{\bef},\\
\int_{\tore} \bw  =0, \int_{\tore} q &= 0, \\
\bw(\bx, 0) = \bw_{0}(\bx)& =\overline{\bv_{0}}, \label{thetallBLM}
\end{align}
with periodic boundary conditions.

\begin{Def} 
Let  ${\bef} \in L^{2}(0,T;\vec{H}^{-1}_{})$ be a divergence free function and  $\bv_0 \in L^{2}_{\sigma }$.  For any $0 \le \theta \le 1$ 	and $ 0 \le N <\infty  $ we say that the couple $(\bw, q) $  is a `` regular" weak solution to (\ref{thetallBM})-(\ref{thetallBLM}) if the following
properties are satisfied:
\begin{align}
\bw &\in \mathcal{C}_{}(0,T;\vec{H}^{\theta}_{\sigma }) \cap
L^2(0,T;\vec{H}^{1+\theta}_{\sigma }),\label{bv12}\\
\bw_{,t}&\in   L^{2}(0,T;\vec{H}^{4 \theta - \frac{3}{2}}_{}),
\label{bvt}\\
q&\in L^{2}(0,T;{H}^{4 \theta - \frac{1}{2}}(\tore))
\label{psp},
\end{align}
the couple $(\bw, q)$ fulfill
\begin{equation}
\begin{split}
\int_0^T \langle \bw_{,t}, \bfi \rangle  -  \langle  \overline{D_{N,\theta} \bw \times \nabla \times D_{N,\theta}\bw}, 
\bfi \rangle  +  \nu \langle 
\nabla \bw, \nabla \bfi  \rangle  +\langle \nabla q,  \bfi \rangle  \; dt\\
 =\int_0^T  \langle \overline{\bef}, \bfi \rangle \; dt
\qquad \textrm{ for all } {\bfi}\in L^{2}(0,T; \vec{H}^{1}).
\end{split}\label{weak1}
\end{equation}
Moreover,
\begin{equation}
\bw(0)= \bw_0.
\label{intiale}
\end{equation}
\end{Def}

Note that the terms $\langle \bw_{,t}, \bfi \rangle$,    $\langle  \overline{D_{N,\theta} \bw \times \nabla \times D_{N,\theta}\bw}, 
\bfi \rangle $ and $\langle \nabla q,  \bfi \rangle$ in (\ref{weak1}) are defined in the sense of duality pairing on $ \vec{H}^{4 \theta - \frac{3}{2}} \times \vec{H}^{\frac{3}{2}-4\theta} $. The term    $\nu \langle 
\nabla \bw, \nabla \bfi  \rangle$  is defined in the sense   of duality pairing on $   \vec{H}^{\theta}\times  \vec{H}^{-\theta}$.
The last term $\langle \overline{\bef}, \bfi \rangle$ is defined in the sense   of duality pairing on $\vec{H}^{1-2\theta} \times \vec{H}^{2\theta-1}$

\begin{Rem}
The notion of ``regular weak solution" is introduced in \cite{bresslilewandowski}. 
Here, we use the name ``regular"  for the weak  solution since  the velocity part of the solution  $\bw$ does not develop a finite time singularity.\\ 
\end{Rem}

   The main result of this section is  to prove the existence of a unique ``regular" weak solution to the RADM model. 
 
 \begin{Theorem} 
Assume  ${\bef} \in L^{2}(0,T;\vec{H}^{-1}_{})$ be a divergence free function and  $\bv_0 \in L^{2}_{\sigma }$.  Let $ \theta \ge \frac{1}{6}$ 	and  let $ 0 \le N <\infty  $ is given and fixed.  Then problem  (\ref{thetallBM})-(\ref{thetallBLM}) has a unique regular
weak solution.
\label{TH1}
\end{Theorem}

\textbf{Proof of Theorem \ref{TH1}.}
The proof of Theorem \ref{TH1} follows the classical scheme. 
For simplicity we restrict ourselves  to the  critical case $\theta = \frac{1}{6}$, 
 and  we drop some indices of $\theta$.  We will consider only the case when $0<N<\infty $ and the case $N=0$ is deduced  by replacing $D_{0,\theta}$ by $I$. Thus  we will write   ``$D_{N}$'' instead of ``$ D_{N,\theta}$'' and ``$\mathbb{A}_{}$'' instead of ``$ \mathbb{A}_{\theta}$''  expecting that no confusion will occur.
 We start by constructing
approximated solutions $(\bw^n, q^n)$ via Galerkin method. Then we seek for a priori
estimates that are uniform with respect to $n$. Next, we passe to the limit in the
equations after having used compactness properties. Finally we show that the solution
we constructed is unique thanks to Gronwall's lemma \cite{RT83}.\\

\textbf{Step 1}(Galerkin approximation).
Consider a sequence $\left\{ \bfi^{r} \right\}_{r=1}^{\infty}$ consisting of $L^2$-orthonormal  and $H^{1}$-orthogonal eigenvectors of the Stokes problem subjected to the
space periodic conditions. 
We note that this sequence forms a hilbertian basis of $L^2$.\\
We set 
\begin{equation}
\begin{split}
\bw^n(t,\bx)=\sum_{{r=1}}^{N}\bc_{r}^n (t) \bfi^{r}(\bx),\\
\quad  \hbox{ and } q^n(t,\bx)=\sum_{{|\vec{k}|=1}}^{n}q_{\vec{k}}^n (t) e^{i \vec{k} \cdot \bx}.
\end{split}
\end{equation}
We look for $(\bw^n(t,\bx), q^n(t,\bx)) $ that are determined through the system of equations 
\begin{equation}
\begin{split}
 \left( \bw_{,t}^n, \bfi^{r} \right)  -  (\overline{D_N\bw^n \times \nabla \times D_N\bw^n}, 
\bfi^{r}) +  \nu(
\nabla \bw^n, \nabla \bfi^{r}  )\; \\
= \langle \bef, \bfi^{r} \rangle \; ,
\quad {r=1},2,...,n,
\end{split}\label{weak1galerkine}
\end{equation}
and 
\begin{equation}
\begin{split}
 \displaystyle \Delta^{} q^n=  \nabla \cdot   \Pi^n \left(  \overline{D_{N} {\bw}^n  \times \nabla \times D_{N} {\bw}^n } \right).
\end{split}\label{pressuregalerkine}
\end{equation}
Where the projector $ \displaystyle \Pi^n $ assign to any Fourier series $\displaystyle \sum_{\bk \in \Z^3\setminus\{0\}} \vec{g}_{\bk} e^{i\bk \cdot\bx} $ the following  series  
 $\displaystyle \sum_{\bk \in \Z^3\setminus\{0\}, |\bk| \le n} \vec{g}_{\bk} e^{i\bk \cdot\bx}. $

Moreover we require that $\bw^n $ satisfies the following initial condition
\begin{equation}
 \label{initial Galerkine}
\bw^n(0,.)= \bw^n_0= \sum_{r=1}^{N}\bc_0^n  \bfi^{r}(\bx),
\end{equation}
and 
\begin{equation}
 \label{initial2 Galerkine}
\bw^n_0 \rightarrow \bw_0  \quad \textrm{ strongly  in }  \vec{H}^{\frac{1}{6}}_{\sigma} \quad \textrm{ when } n \rightarrow \infty .
\end{equation}
 The classical Caratheodory theory \cite{Wa70}  then implies the short-time existence of solutions 
to (\ref{weak1galerkine})-(\ref{pressuregalerkine}).  Next we derive  estimates on $\bc^n$ that is uniform w.r.t. $n$.
These estimates then imply that the  solution of  (\ref{weak1galerkine})-(\ref{pressuregalerkine}) constructed on a short time interval $[0, T^n[ $ exists for all $t \in [0, T]$.\\

\textbf{Step 2} (A priori estimates)
Multiplying the $r$th equation in (\ref{weak1galerkine}) with $\alpha^{\frac{1}{3}}|\bk|^{\frac{1}{3}}\widehat{D_N}{\bc}^n_{r}(t)+\widehat{D_N}{\bc}^n_{r}(t)$, summing over ${r=1},2,...,n$, integrating over time from $0$ to $t$ and using the following identities
 \begin{equation}
 \label{divergencfreebar1}
\left({\bw}_{,t}^n, \mathbb{A}D_N{\bw}^n\right)=
\frac{1}{2}\frac{d}{dt}\|\mathbb{A}^{\frac{1}{2}}D_N^{\frac{1}{2}}{\bw}^n \|_{2}^2 
\end{equation}
\begin{equation}
 \label{divergencfreebar2}
\left(-\Delta{\bw}^n,   \mathbb{A}^{}D_N^{}{\bw}^n\right)=
\|\mathbb{A}^{\frac{1}{2}}D_N^{\frac{1}{2}}{\bw}^n\|_{1,2}^2,
\end{equation}
\begin{equation}
 \label{divergencfreebar3}
\langle \overline{\bef}, \mathbb{A}^{}D_N^{}{\bw}^n \rangle= \langle\mathbb{A}^{\frac{1}{2}}D_N^{\frac{1}{2}} \overline{\bef}, \mathbb{A}^{\frac{1}{2}}D_N^{\frac{1}{2}}{\bw}^n \rangle,
\end{equation}
and
\begin{equation}
 \label{divergencfreebar}
 \begin{array}{lll}
\left(\overline{{D_N\bw^n} \times \nabla \times {D_N\bw^n}},   \mathbb{A}^{}D_N^{}\bw^n\right)&=\left(  D_N\bw^n \times \nabla \times {D_N\bw^n},  D_N^{}\bw^n\right) = 0,
\end{array}
\end{equation}
 leads to the a priori estimate
\begin{equation}
\begin{array}{lllll}
 \label{apriori1}
\displaystyle \frac{1}{2}\| \mathbb{A}^{\frac{1}{2}}D_N^{\frac{1}{2}}{\bw}^n \|_{2}^2 
+ \displaystyle \nu\int_{0}^{t}\| \mathbb{A}^{\frac{1}{2}}D_N^{\frac{1}{2}} {\bw}^n \|_{1,2}^2  \ ds\\
 \quad = \displaystyle \int_{0}^{t} \langle\mathbb{A}^{\frac{1}{2}}D_N^{\frac{1}{2}} \overline{\bef}, \mathbb{A}^{\frac{1}{2}}D_N^{\frac{1}{2}}{\bw}^n \rangle \ ds 
 +\displaystyle \frac{1}{2}\| \mathbb{A}^{\frac{1}{2}}D_N^{\frac{1}{2}}\overline{\bv}^{n}_{0} \|_{2}^2.
\end{array}
\end{equation}
 Using the duality norm combined with Young inequality and inequality (\ref{lemme13})   we conclude from  (\ref{apriori1}) that 
 \begin{equation}
 \label{iciapriori12}
 \begin{split}
\sup_{t \in [0,T^n[}\| \mathbb{A}^{\frac{1}{2}}D_N^{\frac{1}{2}}{\bw}^n \|_{2}^2 + \nu\int_{0}^{t} \| \mathbb{A}^{\frac{1}{2}}D_N^{\frac{1}{2}}{\bw}^n \|_{1,2}^2  \ ds  \le \| {\bv}^{n}_{0} \|_{2}^2 + \frac{1}{\nu}\int_{0}^{T} \| {\bef} \|_{-1,2}^2  \ ds
\end{split}
\end{equation}
that immediately implies that the existence time is independent of $n$ and it is possible to take $T=T^n$.\\ 
We deduce from (\ref{iciapriori12}) and (\ref{lemme13d}) that 
 \begin{equation}
\label{vbar1sansdn}
 {\bw}^n \in L^{\infty}(0,T ; \vec{H}^{\frac{1}{6}}_{\sigma }) \cap L^{2}(0,T ; \vec{H}^{\frac{7}{6}}_{\sigma }).
 \end{equation}
 Thus it follow from (\ref{lemme13a}) that 
 \begin{equation}
\label{icivbar1}
 D_N{\bw}^n \in L^{\infty}(0,T ; \vec{H}^{\frac{1}{6}}_{\sigma }) \cap L^{2}(0,T ; \vec{H}^{\frac{7}{6}}_{\sigma }),
 \end{equation}
 
  From  ({\ref{icivbar1}}) and by using H\"{o}lder inequality combined with Sobolev injection we get 
  \begin{equation}
\label{icivtilde1rev}
{D_N{\bw}^n \times \nabla \times  D_N{\bw}^n} \in L^{2}(0,T ; H^{\frac{-7}{6}}_{}).  
 \end{equation}
 From (\ref{icivtilde1rev}) and (\ref{lemme2.2})  it follows   that  
 \begin{equation}
\label{vtilde1}
\overline{D_N{\bw}^n \times \nabla \times D_N{\bw}^n} \in L^{2}(0,T ; H^{\frac{-5}{6}}_{}).  
 \end{equation}
Consequently from the elliptic theory eqs (\ref{pressuregalerkine}) implies that 
 \begin{equation}
\label{icivbarvbarpressure}
\int_{0}^{T}\|q^n\|_{\frac{1}{6},2}^2 dt < K. 
\end{equation}
From eqs. (\ref{weak1galerkine}), (\ref{vbar1sansdn}), (\ref{vtilde1}) and (\ref{icivbarvbarpressure}) we also obtain that 
 \begin{equation}
\label{icivtemps}
\int_{0}^{T}  \|\bw^n_{,t}\|_{-\frac{5}{6},2}^2  dt < K. 
\end{equation}

\textbf{Step 3} (Limit $n \rightarrow \infty$) It follows from the estimates (\ref{vbar1sansdn})-(\ref{icivtemps}) and the Aubin-Lions compactness lemma
(see \cite{sim87} for example) that there are a  not relabeled  subsequence of $(\bw^n, q^n)$  and a couple $(\bw, q)$ such that
\begin{align}
\bw^n &\rightharpoonup^* \bw &&\textrm{weakly$^*$ in } L^{\infty}
(0,T;\vec{H}^{\frac{1}{6}}_{\sigma }), \label{c122}\\
D_N \bw^n &\rightharpoonup^* D_N \bw &&\textrm{weakly$^*$ in } L^{\infty}
(0,T;\vec{H}^{\frac{1}{6}}_{\sigma }), \label{DNc122}\\
\bw^n &\rightharpoonup \bw &&\textrm{weakly in }
L^2(0,T;\vec{H}^{\frac{7}{6}}_{\sigma }), \label{nc22}\\
D_N\bw^n &\rightharpoonup D_N\bw &&\textrm{weakly in }
L^2(0,T;\vec{H}^{\frac{7}{6}}_{\sigma }), \label{DNnc22}\\
\bw^n_{,t}&\rightharpoonup \bw_{,t} &&\textrm{weakly in } L^{2}
(0,T;\vec{H}^{-\frac{5}{6}}_{}),
\label{nc322}\\
q^n&\rightharpoonup q &&\textrm{weakly in } L^{2}(0,T;H^{\frac{1}{6}}
(\tore)), \label{c32}\\
\bw^n &\rightarrow \bw &&\textrm{strongly in  }
L^2(0,T;\vec{H}^{\frac{1}{6}}_{\sigma }),
\label{c83icil2}\\
D_N\bw^n &\rightarrow D_N\bw &&\textrm{strongly in  }
L^2(0,T;\vec{H}^{\frac{1}{6}}_{\sigma }),
\label{DNc83icil2}
\end{align}
 From (\ref{DNnc22}) and (\ref{DNc83icil2})  it follows   that  
 
 \begin{align}
 \overline{D_N{\bw}^n \times \nabla \times D_N {\bw}^n} &\rightarrow \overline{D_N {\bw} \times \nabla \times D_N\bw} &&\textrm{strongly in  }
L^1(0,T;L^{1}(\tore)^{3}),\label{c82int}
 \end{align}
Finally, since the
sequence $  \overline{D_N{\bw}^n \times \nabla \times D_N{\bw}^n}$  is bounded in $L^2(0,T;H^{-\frac{5}{6}}(\tore)^{3})$, it converges weakly, up to a
subsequence, to some $\chi $ in  $ L^2(0,T;H^{-\frac{5}{6}}(\tore)^{3})$ . The result above and uniqueness of the limit,
allows us to claim that $\chi  =\overline{D_N\bw \times \nabla \times D_N {\bw}} $. Consequently
 \begin{align}
     \overline{D_N{\bw}^n \times \nabla \times D_N{\bw}^n} &\rightharpoonup \overline{D_N\bw \times \nabla \times D_N {\bw}} &&\textrm{weakly in }
 L^2(0,T;H^{-\frac{5}{6}}(\tore)^{3}), \label{c22primen}\    
  \end{align}
The above established convergences are clearly sufficient for taking the limit in (\ref{weak1galerkine})  and for concluding that  the velocity part  $ \bw$  satisfy (\ref{weak1}). 
Moreover, 
from (\ref{nc22}) and (\ref{nc322})
 one we can deduce by a classical argument ( see in \cite{A01})   that 
 \begin{equation}
 \bw \in  \mathcal{C}(0,T;\vec{H}^{\frac{1}{6}}_{\sigma }).
\end{equation}
Furthermore, from  the  strong continuity of $\bw$ with respect to the time with value in $\vec{H}^{\frac{1}{6}}$  we deduce   that $\bw(0)=\bw_0$.\\
Let us mention also that $D_N{\bw}+ \alpha^{^\frac{1}{3}}(-\Delta)^\frac{1}{6}D_N{\bw} \in L^2(0,T;\vec{H}^{\frac{5}{6}}_{\sigma }), $ hence  $\mathbb{A}^{}D_N^{}{\bw} $ is a possible  test function in the weak formulation (\ref{weak1}). Thus $ \mathbb{A}^{\frac{1}{2}}D_N^{\frac{1}{2}}{\bw} $ verifies for all $t \in [0,T]$  the following equality   
\begin{equation}
\begin{array}{lll}
 \label{apriori12leray}
\displaystyle \frac{1}{2}\| \mathbb{A}^{\frac{1}{2}}D_N^{\frac{1}{2}}{\bw} \|_{2}^2 
+ \displaystyle \nu\int_{0}^{t}\| \mathbb{A}^{\frac{1}{2}}D_N^{\frac{1}{2}} {\bw} \|_{1,2}^2  \ ds\\
 \quad = \displaystyle \int_{0}^{t} \langle\mathbb{A}^{\frac{1}{2}}D_N^{\frac{1}{2}} \overline{\bef}, \mathbb{A}^{\frac{1}{2}}D_N^{\frac{1}{2}}{\bw} \rangle \ ds 
 +\displaystyle \frac{1}{2}\| \mathbb{A}^{\frac{1}{2}}D_N^{\frac{1}{2}}\overline{\bv}_{0} \|_{2}^2.
\end{array}
\end{equation}

 \textbf{Step 4} (Uniqueness)
Since the pressure part of the solution is uniquely determined by the velocity part it remain to show the uniqueness to the velocity. 

Next, we will show the continuous dependence of the  solutions on the initial data and in particular the uniqueness.\\
Let $ \theta \ge \frac{1}{6}$  and let $({\bw_1,q_1})$ and $({\bw_2,q_2})$   be any two solutions of (\ref{thetallBM})-(\ref{thetallBLM}) on the interval $[0,T]$, with initial values $\bw_1(0)$ and $\bw_2(0)$. Let us denote by  $\delta \vec{w}_{} =\bw_2-\bw_1$.
 We subtract the equation for $\bw_1$ from the equation for $\bw_2$ and test it with $ \mathbb{A}^{}D_N^{} (\delta \vec{w})$.
 We get using successively the fact that the averaging operator commutes with differentiation under periodic boundary conditions, the norm duality, Young inequality, Lemma \ref{lemme1}:
 \begin{equation}
 \label{notuniformzero}
 \begin{split}
  \displaystyle 
 \frac{1}{2}\frac{d}{dt} \|\mathbb{A}^{\frac{1}{2}}D_N^{\frac{1}{2}}{\delta \vec{w}}\|_{2}^{2} +\nu \|\nabla \mathbb{A}^{\frac{1}{2}}D_N^{\frac{1}{2}}{\delta \vec{w}}_{}\|_{2}^{2} \\
   \le \displaystyle (\overline{D_N(\bw_{2})\times \nabla \times D_N({\bw}_{2})}-\overline{D_N(\bw_{1})\times \nabla \times D_N({\bw}_{1})},   \mathbb{A}^{}D_N^{} (\delta \vec{w}) )
    \\   \le \displaystyle
    ( {D_N\bw_{2}\times \nabla \times  {D_N\bw}_{2}}-{D_N\bw_{1}\times \nabla \times {D_N\bw}_{1}}, {D_N\delta \vec{w}} )\\
\le \displaystyle \frac{C(\alpha)(N+1)^2}{\nu} \|{D_N \delta \vec{w}}\times \nabla \times {D_N \bw}_{1} \|_{-{\frac{7}{6}},2}^{2} + \frac{\nu \alpha^{\frac{1}{3}}}{2} \|{\delta \vec{w}} \|_{\frac{7}{6},2}^2. \end{split}
\end{equation}
 By using  (\ref{lemme13d}) we have 
\begin{equation}
\label{notuniform}
 \begin{split}
  \displaystyle 
   \frac{1}{2}\frac{d}{dt}\alpha^{\frac{1}{3}} \|{\delta \vec{w}}\|_{\frac{1}{6},2}^{2} + \nu \alpha^{\frac{1}{3}} \|\nabla {\delta \vec{w}}_{}\|_{\frac{1}{6}}^{2} \\
\le 
  \frac{1}{2}\frac{d}{dt} \|\mathbb{A}^{\frac{1}{2}}D_N^{\frac{1}{2}}{\delta \vec{w}}\|_{2}^{2} +\nu \|\nabla \mathbb{A}^{\frac{1}{2}}D_N^{\frac{1}{2}}{\delta \vec{w}}_{}\|_{2}^{2} \end{split}
\end{equation}
From (\ref{notuniformzero}) (\ref{notuniform}) and  using (\ref{lemme13a}), H\"{o}lder inequality and Sobolev embedding theorem we get 
\begin{equation}
 \begin{split}
  \displaystyle 
  \frac{d}{dt}\alpha^{\frac{1}{3}}  \|{\delta \vec{w}}\|_{\frac{1}{6},2}^{2} + \nu \alpha^{\frac{1}{3}} \|\nabla {\delta \vec{w}}_{}\|_{\frac{1}{6}}^{2}  \\
 \le \displaystyle \frac{C(\alpha)(N+1)^6}{\nu} \|\delta \vec{w}\|_{\frac{1}{6}}^{2} \| {\bw}_{1} \|^{2}_{\frac{7}{6}}.
\end{split}
\end{equation}

Using Gronwall's inequality, since $\| {\bw}_{1} \|^{2}_{\frac{7}{6}} \in L^1([0,T])$, we conclude 
the continuous dependence of the solutions on the initial data in the $L^{\infty}([0,T],\vec{H}^{\frac{1}{6}}_{\sigma })$  norm. In particular, if ${\delta \vec{w}}^{}_{0}=0$ then ${\delta \vec{w}}=0$ and the solutions are unique for all $t \in [0,T] .$  Since $T>0$ is arbitrary this solution may be uniquely extended for all time.\\ 
This finishes the proof of Theorem \ref{TH1}.\\

\section{Limit when $N \rightarrow \infty $  }
\hskip 0.5cm 
Let $(\vec{w}_N, q_N ) $   be the solution of the RADM model constructed above.  In this section, we  take the limit $N \rightarrow \infty $ and look for the
 optimal value of $\theta$ 
    in order to show that  $(\vec{w}_N, q_N ) $  converges,  up to subsequences, to a solution of the mean 
rotational  Navier-Stokes equations. 

 \begin{Theorem}
\label{1deuxieme}
 Let $(\bw_{N}, q_{N})$ be the solution of (\ref{thetallBM})-(\ref{thetallBLM}) for a fixed $\alpha>0$, $ N>0$ and  $  \frac{1}{6} \le  \theta <1$. 
 There is a subsequence $N_j \rightarrow  \infty $ as  $j \rightarrow \infty $  such that:
$(\bw_{N_j}, q_{N_j}) \rightarrow (\bw, q)$ as $j \rightarrow \infty $
where
$ (\bw, q)  \in L^{\infty}
([0,T];\vec{H}^{\theta})\cap L^2([0,T];\vec{H}^{1+\theta}_{}) \times L^2([0,T];H^{-\frac{1}{2}+2\theta}(\tore)) $ is a distributional solution of  the following system with periodic boundary conditions  
\begin{align}
\nabla \cdot  \bw_{} &=0, \label{AthetallBM}\\
\bw_{,t} -  \overline{\mathbb{A}_{\theta} \bw \times \nabla \times \mathbb{A}_{\theta}\bw} -  \nu \Delta \bw + \nabla q
 &=  \overline{\bef},\\
\int_{\tore} \bw  =0, \int_{\tore} q &= 0, \\
\bw(\bx, 0) = \bw_{0}(\bx)& =\overline{\bv_{0}}. \label{AthetallBLM}
\end{align}
 The sequence  $ \bw_{N_j}$  converges strongly
to $\bw$ in the space $L^2([0,T];\vec{H}^{s})  $ for all  $  s< 1+\theta, $   while the sequence
 $  q_{N_j} $  converges weakly to $q$ in the space $ L^2([0,T];H^{-\frac{1}{2}+2\theta}(\tore)).$
\end{Theorem}
\begin{Rem}
(1) The restriction in the case of the ADM model  in \cite{bresslilewandowski} is $\theta > \frac{3}{4}$.
Our result is sharper for the RADM model and it should be clear to the experts that the same approach used here can be used to the ADM model to get similar results.\\
(2) For $\theta \ge 1$, one can uses the same method developed in \cite{bresslilewandowski} to show the convergence of RADM model to the    mean rotational  Navier-Stokes equations. 
\end{Rem}
In order to prove Theorem \ref{1deuxieme} we need to
reconstruct a uniform estimates for $\bw_{N}$. Since some of the a priori  estimates used  in the proof of Theorem \ref{TH1}    dependent   on $N$  we can not use them. Instead we  show in the proof of Theorem \ref{1deuxieme}  that $D_{N,\theta}\bw_{N}$ belongs to $ L^{\infty}([0,T];\vec{L}^{2})\cap L^2([0,T];\vec{H}^{1}_{\sigma })$ uniformly with respect to $N$.\\

Before proving Theorem \ref{1deuxieme}, we first record the following  Lemma.
%
\begin{Lem}
\label{fourierdiscret}
Let $ \theta >0 $ and assume that $\bv  \in L^{2}([0,T], \vec{H}^{2 \theta} )$.   Then   
\begin{align} D_{N,\theta}  \bv  
  &\rightarrow  \mathbb{A}_{\theta} \bv  \textrm{ strongly in }
L^2(0,T;\vec{L}^2), \end{align}
and  there exist a constant C  independent from $N$  such that \begin{align} \| D_{N,\theta}  \bv \|_{2}  \le \| \mathbb{A}_{\theta} \bv\|_{2 }  \le C(\alpha^{\theta}) \|  \bv\|_{2\theta,2 }.    \end{align}
\end{Lem}
\textbf{Proof.} The first part of this lemma is given page 18 in \cite{bresslilewandowski}. The second part  is a direct consequence from the propriety  (\ref{haniali}) of the operator $D_{N,\theta}$, the fact that $$ \| \mathbb{A}_{\theta} \bv\|_{2 }^2= \sum_{{\vec{k} \in {\mathcal T}_{3}}}^{} \left(1+\alpha^{2\theta}|{\bk}|^{2\theta}\right)^2|{\bc}_{\bk}|^2=\| \bv\|_{2 }^2 + 2\alpha^{2\theta} \| \bv\|_{\theta,2 }^2 + \alpha^{4\theta}\| \bv\|_{2\theta,2 }^2,$$  and Poincaré inequality.

\textbf{Proof of Theorem \ref{1deuxieme}.}
  The proof of Theorem \ref{1deuxieme} follows the lines of the proof of the Theorem 4.1 in \cite{bresslilewandowski},
  that we have to modify in order to treat the cases when $\theta \le \frac{3}{4}$. 
Thus  we will use  \cite{bresslilewandowski} as a
reference and only point out the differences between their  proof of convergence to the mean 
 Navier-Stokes equations and the proof of convergence  in the present study.  
First, we need to find estimates
that are independent from $N$. \\

\textbf{Step 1} (Uniform estimates with respect to N)
Following  the proof of Theorem \ref{TH1}  we obtain  
that the solution of (\ref{thetallBM})-(\ref{thetallBLM})  satisfies 
\begin{equation}
\begin{array}{lll}
 \label{undemiapriori12leray}
\displaystyle \frac{1}{2}\| \mathbb{A}_{\theta}^{\frac{1}{2}}D_{N,\theta}^{\frac{1}{2}}{\bw_N} \|_{2}^2 
+ \displaystyle \nu\int_{0}^{t}\| \mathbb{A}_{\theta}^{\frac{1}{2}}D_{N,\theta}^{\frac{1}{2}} {\bw_N} \|_{1,2}^2  \ ds\\
 \quad \le \displaystyle \int_{0}^{t} \langle\mathbb{A}_{\theta}^{\frac{1}{2}}D_{N,\theta}^{\frac{1}{2}} \overline{\bef}, \mathbb{A}_{\theta}^{\frac{1}{2}}D_{N,\theta}^{\frac{1}{2}}{\bw_N} \rangle \ ds 
 +\displaystyle \frac{1}{2}\| \mathbb{A}_{\theta}^{\frac{1}{2}}D_{N,\theta}^{\frac{1}{2}}\overline{\bv}_{0} \|_{2}^2.
\end{array}
\end{equation}
Thus 

 \begin{equation}
 \label{apriori12}
 \begin{split}
\sup_{t \in [0,T]}\| \mathbb{A}_{\theta}^{\frac{1}{2}}D_{N,\theta}^{\frac{1}{2}}{\bw}_N \|_{2}^2 + \nu\int_{0}^{t} \| \mathbb{A}_{\theta}^{\frac{1}{2}}D_{N,\theta}^{\frac{1}{2}}{\bw}_N \|_{1,2}^2  \ ds  \le \| {\bv}^{}_{0} \|_{2}^2 + \frac{1}{\nu}\int_{0}^{T} \| {\bef} \|_{-1,2}^2  \ ds
\end{split}
\end{equation}
 consequently we can bound the right hand side by a constant C which  is independent from $N$.

We deduce from (\ref{apriori12}) and Lemma \ref{lemme1} that 
\begin{equation}
\label{vbar1}
 {D_{N,\theta}\bw_N} \in L^{\infty}(0,T ; \vec{L}^{2}_{\sigma }) \cap L^{2}(0,T ; \vec{H}^{1}_{\sigma }), \textrm{ uniformly with respect to }  N,
 \end{equation}
and 
\begin{equation}
\label{vbar2}
 {\bw_N} \in L^{\infty}(0,T ; \vec{H}^{\theta}_{\sigma }) \cap L^{2}(0,T ; \vec{H}^{1+\theta}_{\sigma }), \textrm{ uniformly with respect to } N.
 \end{equation}
 
 We observe from (\ref{vbar1}) that
 \begin{equation}
\label{DNvbarsansbar}
{D_{N,\theta} ({\bw_N}) \times \nabla \times  D_{N,\theta} ({\bw_N})} \in  L^{2}(0,T ;{H}^{-\frac{3}{2}}(\tore)^{3 } ).
\end{equation}
Thus from (\ref{DNvbarsansbar}) and  (\ref{VPYRA}) we obtain  
 \begin{equation}
\label{vbarvbar}
\overline{{D_{N,\theta} ({\bw_N}) \times \nabla \times  D_{N,\theta} ({\bw_N})}} \in  L^{2}(0,T ;{H}^{-\frac{3}{2} + 2\theta }(\tore)^{3} ).
\end{equation}
For the pressure term $q_N$, we deduce  that it verifies the following equation 
\begin{equation}
\begin{split}
\displaystyle \Delta^{} q_N =  \nabla \cdot    \left(  \overline{D_{N,\theta} {\bw}_N  \times \nabla \times D_{N,\theta} {\bw}_N }       \right).
\end{split}\label{pressure}
\end{equation}
consequently  the classical elliptic theory combined with (\ref{vbarvbar}) implies that 
 \begin{equation}
\label{vbarvbarpressure}
\int_{0}^{T}\|q_{N}\|_{-\frac{1}{2} + 2\theta,2}^2 dt < K, \textrm { uniformly with respect to  } N 
\end{equation}
From (\ref{thetallBM})-(\ref{thetallBLM})  and (\ref{vbarvbar}) we also obtain that 
 \begin{equation}
\label{vtemps}
\int_{0}^{T}  \|\partial_t\bw_{N}\|_{-\frac{3}{2} + 2\theta,2}^2  dt < K, \textrm{ uniformly with respect to } N. 
\end{equation}

\textbf{Step 2} (Passing to the limit $ N \rightarrow \infty$) The central issues is how to take the limit in the nonlinear term $ {D_{N,\theta} ({\bw_N}) \times \nabla \times  D_{N,\theta}({\bw_N})}$.

From the  Aubin-Lions compactness Lemma  (the same arguments as in section 3) we can find a  subsequence $(\vec{w}_{N_j}, D_{N_j, \theta}\vec{w}_{N_j}, q_{N_j})$  and $(\vec{w}_{}, \vec{z},q_{})$ such that when   $N_j \rightarrow \infty $  we have: 

\begin{align}
\bw_N &\rightharpoonup^* \bw &&\textrm{weakly$^*$ in } L^{\infty}
(0,T;\vec{H}^{\theta}_{\sigma }), \label{limitNc122}\\
\bw_N &\rightharpoonup \bw &&\textrm{weakly in }
L^2(0,T;\vec{H}^{1+\theta}_{\sigma }), \label{limitNc22}\\
D_{N,\theta}\bw_N &\rightharpoonup \vec{z} &&\textrm{weakly in }
L^2(0,T;\vec{H}^{1}_{\sigma }), \label{limitNDNc22}\\
\partial_t \bw_{N}&\rightharpoonup \partial_t \bw &&\textrm{weakly in } L^{2}
(0,T;\vec{H}^{-\frac{3}{2} + 2\theta}_{}),
\label{limitNc322}\\
q_N&\rightharpoonup q &&\textrm{weakly in } L^{2}(0,T;H^{-\frac{1}{2}+2\theta}
(\tore)), \label{limitNc32}\\
\bw_N &\rightarrow \bw &&\textrm{strongly in  }
L^2(0,T;\vec{H}^{s}_{\sigma }) \textrm{ for all } s < 1+\theta.
\label{limitNc83icil2}
\end{align}

 The goal is
to prove 
 \begin{align}
 D_{N,\theta}{\bw}_N \times \nabla \times D_{N,\theta} {\bw}_N &\rightarrow  \mathbb{A}_{\theta}\bw \times \nabla \times \mathbb{A}_{\theta}\bw &&\textrm{strongly in  }
L^1(0,T;L^{1}(\tore)^{3}),\label{DNc82int}
 \end{align}
Thus 
 it remain to show that 
\begin{align}
 D_{N,\theta}{\bw}_N &\rightarrow  \mathbb{A}_{\theta} {\bw}  &&\textrm{strongly in  }
L^2(0,T;\vec{L}^{2}).\label{hhdac82int}
 \end{align}
 In order to show (\ref{hhdac82int}), we compute directly the difference between  $ D_{N,\theta}{\bw}_N $ and  $\mathbb{A}_{\theta} {\bw}$ as follows 
\begin{equation}  
 \begin{split}
\| D_{N,\theta}{\bw}_N - \mathbb{A}_{\theta}{\bw}\|_{L^2(0,T;\vec{L}^{2})} \le \| D_{N,\theta}{\bw}_N - D_{N,\theta}{\bw}\|_{L^2(0,T;\vec{L}^{2})} + \| D_{N,\theta}{\bw}-  \mathbb{A}_{\theta} {\bw}  \|_{L^2(0,T;\vec{L}^{2})}\\
\le \| \mathbb{A}_{\theta} \left( {\bw}_N - {\bw} \right) \|_{L^2(0,T;\vec{L}^{2})} + \| D_{N,\theta}{\bw}- \mathbb{A}_{\theta}{\bw}\|_{L^2(0,T;\vec{L}^{2})}\\
 \le C\| {\bw}_N -{\bw}\|_{L^2(0,T;\vec{H}^{2\theta})} + \| D_{N,\theta}{\bw}- \mathbb{A}_{\theta}{\bw}\|_{L^2(0,T;\vec{L}^{2})} 
 \end{split}
 \end{equation}
 Where we have used the second part of lemma \ref{fourierdiscret}. 
 Hence, using that $ \theta <1$, (\ref{limitNc22})  and the first part of lemma  (\ref{fourierdiscret}) we deduce (\ref{hhdac82int}). 
 Finally the result above, (\ref{hhdac82int}), combined with (\ref{limitNDNc22})  and the  uniqueness of the limit, allows us to deduce  

 \begin{align}
     \nabla \times D_{N,\theta}{\bw}_N &\rightharpoonup  \nabla \times \mathbb{A}_{\theta} {\bw} &&\textrm{weakly in }
 L^2(0,T;\vec{L}^{2}). \label{c22prime0}\    
  \end{align}

Consequently we get (\ref{DNc82int}) and  by using (\ref{DNvbarsansbar}) and  (\ref{lemme2.2}) we obtain 

\begin{align}
     \overline{D_{N,\theta}{\bw}_N \times \nabla \times D_{N,\theta}{\bw}_N} &\rightharpoonup \overline{\mathbb{A}_{\theta} \bw \times \nabla \times \mathbb{A}_{\theta} {\bw}} &&\textrm{weakly in }
 L^2(0,T;H^{-\frac{3}{2} + 2\theta}(\tore)^{3}). \label{c22prime}\    
  \end{align}

These convergence results allow us  to prove in the same way as in \cite{bresslilewandowski}  that $ (\bw,q)$ is a distributional  solution to the mean rotational Navier-Stokes equations (\ref{AthetallBM})-(\ref{AthetallBLM}), and  $ (\mathbb{A}_{\theta} \bw,\mathbb{A}_{\theta} q)$  is a  distributional solution to the rotational  Navier-Stokes equations  (\ref{nsBM})-(\ref{nsBLM}), so we will not repeat it.





\end{document}